\documentclass[aip,reprint]{revtex4-1}
\usepackage{color,graphicx,amsmath,amsfonts,amssymb,amsxtra,setspace,xspace}
\usepackage[colorlinks=true]{hyperref}
\hypersetup{urlcolor=blue, citecolor=red, linkcolor=blue}

\newtheorem{thm}{Theorem}[section]

\newtheorem{lem}[thm]{Lemma}
\newtheorem{cor}[thm]{Corollary}

\newcommand{\be}[1]{\begin{equation}\label{#1}}
\newcommand{\ee}{\end{equation}}
\renewcommand{\(}{\left(}
\renewcommand{\)}{\right)}
\newcommand{\finprf}{\unskip\null\hfill$\;\square$\vskip 0.3cm}

\newcommand{\R}{{\mathbb R}}
\renewcommand{\S}{{\mathbb S}}
\newcommand{\N}{{\mathbb N}}

\newcommand{\nrmS}[2]{\|{#1}\|_{\mathrm L^{#2}(\S^1)}}
\newcommand{\intpi}[1]{\int_{-\pi}^\pi{#1}\,\mathrm d\sigma}

\begin{document}
\title{Magnetic rings}

\author{Jean Dolbeault}
\author{Maria J.~Esteban}
\affiliation{CEREMADE (CNRS UMR n$^\circ$ 7534), PSL research university, Universit\'e Paris-Dauphine, Place de Lattre de Tassigny, 75775 Paris 16, France }
\email{dolbeaul@ceremade.dauphine.fr}
\email{esteban@ceremade.dauphine.fr}

\author{Ari Laptev}
\affiliation{Department of Mathematics, Imperial College London, Huxley Building, 180 Queen's Gate, London SW7 2AZ, UK }
\altaffiliation[Also at ]{Department of Mathematics, Siberian Federal University, Russia }
\email{a.laptev@imperial.ac.uk}

\author{Michael Loss}
\affiliation{School of Mathematics, Skiles Building, Georgia Institute of Technology, Atlanta GA 30332-0160, USA }
\email{loss@math.gatech.edu}

\date{\today}

\begin{abstract}
We study functional and spectral properties of perturbations of the operator $-(\partial_s+i\,a)^2$ in $\mathrm L^2(\S^1)$. This operator appears when considering the restriction to the unit circle of a two dimensional Schr\"odinger operator with the Bohm-Aharonov vector potential. We prove a Hardy-type inequality on $\R^2$ and, on $\S^1$, a sharp interpolation inequality and a sharp Keller-Lieb-Thirring inequality.
\end{abstract}

%
%


\maketitle


\section{Introduction}\label{Sec:intro}

On the two-dimensional Euclidean space $\R^2$, let us introduce the polar coordinates $(r,\vartheta)\in[0,+\infty)\times\S^1$ of $\mathbf x\in\R^2$ and consider a magnetic potential $\mathbf a$ in a transversal gauge, or Poincar\'e gauge~\cite{MR1219537}, so that $(\mathbf a,\mathbf e_r)=0$ and $(\mathbf a,\mathbf e_\vartheta)=a_\vartheta(r,\vartheta)$, where $(\mathbf e_r,\mathbf e_\vartheta)$ is the oriented orthogonal basis associated with the polar coordinates such that, for any $\mathbf x\in\R^2\setminus\{\mathbf 0\}$, $\mathbf e_r=\mathbf x/r$, $r=|\mathbf x|$. With this notation, the energy $\int_{\R^2}|(i\,\nabla+\mathbf a)\,\Psi|^2\,\mathrm \mathrm d\mathbf x$ corresponding to the magnetic Schr\"odinger operator $-\Delta_{\mathbf a}$ can be rewritten as
\[
\int_0^{+\infty}\int_{-\pi}^\pi\(\,|\partial_r\Psi|^2+\frac1{r^2}\,|\,\partial_\vartheta\Psi+i\,r\,a_\vartheta\,\Psi|^2\)r\,\mathrm d\vartheta\,\mathrm dr\,.
\]
One of the main motivations is the study of \emph{Bohm-Aharonov magnetic fields}~\cite{MR2129228,LaWe1999} with $a_\vartheta(r,\vartheta)=a/r$ for some constant $a\in\R$. We recall that Stokes' formula applied to the magnetic field $b=\mathrm{curl}\,\mathbf a$ shows that the \emph{magnetic flux} is given by
\[
\int_{|\mathbf x|<r}b\,\mathrm \mathrm d\mathbf x=\frac1{2\pi}\int_{-\pi}^\pi a_\vartheta(r,\vartheta)\,r\,\mathrm d\vartheta=a\,.
\]
The main result concerning Bohm-Aharonov magnetic fields is, for an arbitrary non-negative function $\varphi$ in $\mathrm L^q(\S^1)$, $q\in(1,+\infty)$, the Hardy-type inequality
\be{Hardy1}
\int_{\R^2}|(i\,\nabla+\mathbf a)\,\Psi|^2\,\mathrm \mathrm d\mathbf x\ge\tau\int_{\R^2}\frac{\varphi(\mathbf x/|\mathbf x|)}{|\mathbf x|^2}\,|\Psi|^2\,\mathrm d\mathbf x
\ee
which holds for some constant $\tau$ depending on $\nrmS\varphi q$. A precise statement will be given in Corollary~\ref{Hardy}.

The proof relies on a method~\cite{HoffLa2015} developed recently and uses a \emph{Keller-Lieb-Thirring inequality} for the first eigenvalue of a magnetic Schr\"odinger operator on a \emph{magnetic ring} (see Corollary~\ref{Lambda1}). This spectral estimate is equivalent to \emph{sharp interpolation inequalities} for a magnetic Laplacian on the circle and has been inspired by a series of previous papers~\cite{DoEsLa2013,MR3090125,DoEsLaLo2016} on interpolation inequalities and their spectral counterparts. Let us mention that some semiclassical properties of the spectrum of magnetic rings were recently studied including an electric potential that admits a double symmetric well~\cite{BNoFrRa} (also see earlier references therein). Our results are not limited to the semi-classical regime.

\section{Main results}\label{Sec:main}

On $(-\pi,\pi]\approx\S^1$, let us consider the uniform probability measure $\mathrm d\sigma=\mathrm ds/(2\pi)$ and denote by $\nrmS\psi p$ the corresponding $\mathrm L^p$ norm, for any $p\ge1$. Assume that $a:\R\to\R$ is a $2\pi$-periodic function such that its restriction to $(-\pi,\pi]\approx\S^1$ is in $\mathrm L^1(\S^1)$ and define the subspace
\[
X_a:=\left\{\psi\in C_{\mathrm{per}}(\R)\,:\,\psi'+i\,a\,\psi\in\mathrm L^2(\S^1)\right\}
\]
of the space $C_{\mathrm{per}}(\R)$ of the continuous $2\pi$-periodic functions on $\R$. The change of function
\[
\psi(s)\mapsto e^{i\int_{-\pi}^s\(a(s)-\bar a\)\mathrm d\sigma}\,\psi(s)\,,
\]
where $\bar a:=\int_{-\pi}^\pi a(s)\,\mathrm d\sigma$ is the \emph{magnetic flux}, reduces the problem to the case of a constant: in the sequel of this paper we shall always assume that
\[
\mbox{\emph{$a$ is a constant function}.}
\]
Replacing $\psi$ by $s\mapsto e^{iks}\,\psi(s)$ for any $k\in\mathbb Z$ shows that $\mu_{a,p}(\alpha)=\mu_{k+a,p}(\alpha)$ so that we can restrict the problem to $a\in[0,1]$. By considering $\chi(s)=e^{-is}\,\overline{\psi(s)}$, we find
\[
|\psi'+i\,a\,\psi|^2=|\chi'+i\,(1-a)\,\chi|^2=\big|\overline\psi'-i\,a\,\overline\psi\big|^2\,,
\]
and thus $\mu_{a,p}(\alpha)=\mu_{1-a,p}(\alpha)$: it is thus enough to consider the case $a\in[0,1/2]$. 

Using a Fourier series $\psi(s)=\sum_{k\in\mathbb Z}\psi_k\,e^{iks}$, we obtain that
\[
\nrmS{\psi'+i\,a\,\psi}2^2=\sum_{k\in\mathbb Z}\,(a+k)^2\,|\psi_k|^2\ge a^2\,\nrmS\psi2^2\,,
\]
so that $\psi\mapsto\nrmS{\psi'+i\,a\,\psi}2^2+\alpha\,\nrmS\psi2^2$ is coercive for any $\alpha>-\,a^2$. Moreover, the optimal constant $\mu_{a,p}(\alpha)$ in the interpolation inequality
\be{Ineq:Interp}
\nrmS{\psi'+i\,a\,\psi}2^2+\alpha\,\nrmS\psi2^2\ge\mu_{a,p}(\alpha)\,\nrmS\psi p^2
\ee
written for any $\psi\in X_a$ is an increasing concave function of $\alpha>-\,a^2$ characterized by 
\be{minpsi}
\mu_{a,p}(\alpha):=\inf_{\psi\in X_a\setminus\{0\}}\frac{\int_{-\pi}^\pi\(|\psi'+i\,a\,\psi|^2+\alpha\,|\psi|^2\) \mathrm d\sigma}{\nrmS\psi p^2}
\ee
and~\cite{DoEsLaLo2016} $\lim_{\alpha\to-\,a^2}\mu_{a,p}(\alpha)=0$. The inequality~\eqref{Ineq:Interp} is known if either $p=+\infty$~\cite{MR1330606,Ilyin-Laptev-Loss-Zelik2016} or $p=-\,2$~\cite{Exner-Harrell-Loss98} and the expression of an optimal function was given as a series~\cite{MR1330606} for any $\alpha>-\,a^2$ when $p=+\infty$. Our first result is the extension of this \emph{interpolation} result to the case $p\in(2,+\infty)$.
\begin{thm}\label{main} For any $p>2$, $a\in\R$, and $\alpha>-\,a^2$, the infimum in~\eqref{minpsi} is achieved and
\begin{enumerate}
\item[(i)] if $a\in[0,1/2]$ and $a^2\,(p+2)+\alpha\,(p-2)\le1$, then $\mu_{a,p}(\alpha)=a^2+\alpha$ and equality in~\eqref{Ineq:Interp} is achieved only by the constant functions,
\item[(ii)] if $a\in[0,1/2]$ and $a^2\,(p+2)+\alpha\,(p-2)>1$, then $\mu_{a,p}(\alpha)<a^2+\alpha$ and equality in~\eqref{Ineq:Interp} is not achieved by the constant functions.
\end{enumerate}
Moreover, for any $\alpha>-\,a^2$, $a\mapsto\mu_{a,p}(\alpha)$ is monotone increasing on $(0,1/2)$.
\end{thm}
More can be said on $\mu_{a,p}(\alpha)$: see Theorem~\ref{main2}. The region $a^2\,(p+2)+\alpha\,(p-2)<1$ is exactly the set where the constant functions are linearly stable critical points. See Figs.~\ref{DoEsLaLo-F1} and~\ref{DoEsLaLo-F2}.

\medskip With the results of Theorem~\ref{main} in hand, we study some spectral properties of the magnetic Schr\"odinger operator $H_a-\varphi$ on the unit circle $\S^1\approx(-\pi,\pi]\ni s$ where~$\varphi$ is a potential and $H_a$ is the magnetic Laplacian given by
\[
H_a\psi(s)=-\(\frac d{ds}+i\,a\)^2\psi(s)\,.
\]
The presence of a non-trivial magnetic field $a$ in $H_a$ ``lifts'' the spectrum up and the final result substantially depends on its value. Note that Lieb-Thirring inequalities with magnetic field~\cite{Ilyin-Laptev-Loss-Zelik2016}, in particular, imply an inequality for the first eigenvalue. However, it is not known if the constant is sharp. A somewhat similar result where the lifting of the spectrum is provided by a constant magnetic field was proved with different methods~\cite{DoEsLaLo2016}.

The first spectral consequence of Theorem~\ref{main} is a \emph{Keller-Lieb-Thirring inequality} for the first eigenvalue $\lambda_1(H_a-\varphi)$ of the Schr\"odinger operator $H_a-\varphi$. The function $\alpha\mapsto\mu_{a,p}(\alpha)$ is monotone increasing, concave, and therefore has an inverse, denoted by $\alpha_{a,p}:\R^+\to(-a^2,+\infty)$, which is monotone increasing, and convex.
\begin{cor}\label{Lambda1} Let $p>2$, $a\in\left[0,1/2\right]$, $q=p/(p-2)$ and assume that $\varphi$ is a non-negative function in $\mathrm L^q(\S^1)$. Then
\be{lambda1}
\lambda_1(H_a-\varphi)\ge-\,\alpha_{a,p}\(\nrmS\varphi q\).
\ee
If $4\,a^2+\mu\,(p-2)\le1$, then $\alpha_{a,p}(\mu)=\mu-a^2$; if $4\,a^2+\mu\,(p-2)>1$, then $\alpha_{a,p}(\mu)>\mu-a^2$.

These estimates are optimal in the sense that there exists a non-negative function $\varphi$ such that $\lambda_1(H_a-\varphi)=-\,\alpha_{a,p}\(\nrmS\varphi q\)$. If $4\,a^2+\mu\,(p-2)\le1$, then the equality in~\eqref{lambda1} is achieved by constant potentials.\end{cor}

The second application of Theorem~\ref{main} is related to a \emph{Hardy inequality} in $\R^2$. Let us consider the \emph{Bohm-Aharonov vector potential}
\[
\mathbf a(\mathbf x)=a\(\frac{x_2}{|\mathbf x|^2},\frac{-x_1}{|\mathbf x|^2}\),\quad\mathbf x=(x_1,x_2)\in\R^2\,,\quad a\in\R\,.
\]
and recall the inequality~\cite{LaWe1999}
\be{LaWe}
\int_{\R^2}|(i\,\nabla+\mathbf a)\,\Psi|^2\,\mathrm d\mathbf x\ge\min_{k\in\mathbb Z}\,(a-k)^2\int_{\R^2}\frac{|\Psi|^2}{|\mathbf x|^2}\,\mathrm d\mathbf x\,.
\ee
Using interpolation inequalities~\cite{DoEsLa2013}, the following version~\cite{HoffLa2015} of Hardy's inequality in the case $d\ge3$ was proved:
\[
\int_{\R^d}|\nabla\Psi|^2\,\mathrm d\mathbf x\ge\tau\int_{\R^d}\frac{\varphi(\mathbf x/|\mathbf x|)}{|\mathbf x|^2}\,|\Psi|^2\,\mathrm d\mathbf x\,,
\]
where the constant $\tau$ depends on the value of $\|\varphi\|_{\mathrm L^q(\S^{d-1})}$. Using similar arguments we are now able to prove the following result.
\begin{cor}\label{Hardy} Let $p>2$, $a\in\left[0,1/2\right]$, $q=p/(p-2)$ and assume that $\varphi$ is a non-negative function in $\mathrm L^q(\S^1)$. Then Inequality~\eqref{Hardy1} holds with $\tau>0$ being the unique solution of the equation
\[
\alpha_{a,p}(\tau\,\nrmS\varphi q)=0\,.
\]
Moreover, $\tau=a^2/\nrmS\varphi q$ if $4\,a^2+\nrmS\varphi q\,(p-2)\le1$.\end{cor}
Notice that for any $a\in(0,1/2)$, by taking $\varphi$ constant, small enough in order that $4\,a^2+\nrmS\varphi q\,(p-2)\le1$, we recover the inequality
\[
\int_{\R^2}|(i\,\nabla+\mathbf a)\,\Psi|^2\,\mathrm d\mathbf x \ge a^2\int_{\R^2} \frac{|\Psi|^2}{|\mathbf x|^2}\,\mathrm d\mathbf x\,,
\]
which is a equivalent to~\eqref{LaWe}. The case $a=1/2$ is obtained by a limiting procedure and for arbitrary values of $a\in\R$, we refer to the observations of Section~\ref{proof main}.

\section{Proof of Theorem~\ref{main} and further results}\label{proof main}

\begin{lem}\label{Lem:Existence} For all $a\in\R$, $p\in(2,\infty)$ and $\alpha\ge -a^2$, equality in~\eqref{Ineq:Interp} is achieved by at least one function in $X_a$.\end{lem}
Indeed, by the diamagnetic inequality
\[
\big|\,|\psi|'\big|\le|\psi'+i\,a\,\psi|\quad\mbox{a.e.}\,,
\]
which holds for any $\psi\in X_a$, we infer that any minimizing sequence $\{\psi_n\}$ for~\eqref{minpsi} can be taken bounded in $\mathrm H^1(\S^1)$. By the compact Sobolev embeddings, this sequence is relatively compact in $\mathrm L^p(\S^1)$ and in $C(\S^1)$. The maps $\psi\mapsto\intpi{|\psi|^2}$ and $\psi\mapsto\int_{-\pi}^\pi|\psi'+i\,a\,\psi|^2\,\mathrm d\sigma$ are lower semicontinuous by Fatou's lemma, which proves the claim.\finprf

The minimization problem~\eqref{minpsi} has several reformulations, that have already been used in the case $\alpha=0$~\cite{2017arXiv171208829N}.
\\[4pt] 
1) Any solution $\psi\in X_a$ of the minimization problem~\eqref{minpsi} satisfies the Euler-Lagrange equation
\[
(H_a+\alpha)\,\psi=|\psi|^{p-2}\,\psi
\]
up to a multiplication by a constant. We observe that $v(s)=\psi(s)\,e^{ias}$ satisfies the condition
\be{Cdt:Per}
v(s+2\pi)=e^{2i\pi a}\,v(s)\quad\forall\,s\in\R\,,
\ee
and we can reformulate~\eqref{minpsi} as
\[
\mu_{a,p}(\alpha)=\min_{v\in Y_a\setminus\{0\}}\mathsf Q_{p,\alpha}[v]
\]
where $Y_a:=\{v\in C(\R)\,:\,v'\in\mathrm L^2(\S^1)\,,\;\mbox{\eqref{Cdt:Per} holds}\}$ and
\[
\mathsf Q_{p,\alpha}[v]:=\frac{\nrmS{v'}2^2+\alpha\,\nrmS v2^2}{\nrmS vp^2}\,.
\]
2) With $v=u\,e^{i\phi}$ written in polar form, the boundary condition becomes
\be{bcuphi}
u(\pi)=u(-\pi)\,,\quad\phi(\pi)=2\pi\,(a+k)+\phi(-\pi)
\ee
for some $k\in\mathbb Z$, and $\nrmS{v'}2^2=\nrmS{u'}2^2+\nrmS{u\,\phi'}2^2$. We can reformulate~\eqref{minpsi} as
\[
\mu_{a,p}(\alpha)=\kern-6pt\min_{(u,\phi)\in Z_a\setminus\{0\}}\kern-10pt\frac{\nrmS{u'}2^2+\nrmS{u\,\phi'}2^2+\alpha\,\nrmS u2^2}{\nrmS up^2}
\]
where\\
$Z_a:=\{(u,\phi)\in C(\R)^2\,:\,u',\,u\,\phi'\in\mathrm L^2(\S^1)\,,\;\mbox{\eqref{bcuphi} holds}\}$.
\\[4pt]
3) The third reformulation of~\eqref{minpsi} relies on the Euler-Lagrange equations
\[
-\,u''+|\phi'|^2\,u+\alpha\,u=|u|^{p-2}\,u\quad\mbox{and}\quad(\phi'\,u^2)'=0\,.
\]
Integrating the second equation, and \emph{assuming that $u$ never vanishes}, we find a constant $L$ such that $\phi'=L/u^2$. Taking~\eqref{bcuphi} into account, we deduce from
\[
L\int_{-\pi}^\pi\frac{\mathrm ds}{u^2}=\int_{-\pi}^\pi\phi'\,\mathrm ds=2\pi\,(a+k)
\]
that
\[
\nrmS{u\,\phi'}2^2=L^2\int_{-\pi}^\pi\frac{\mathrm d\sigma}{u^2}=\frac{(a+k)^2}{\nrmS{u^{-1}}2^2}\,.
\]
Hence
\[
\phi(s)-\phi(0)=\frac{a+k}{\nrmS{u^{-1}}2^2}\int_{-\pi}^s\frac{\mathrm ds}{u^2}\,.
\]
Let us define
\[
\mathcal Q_{a,p,\alpha}[u]:=\frac{\nrmS{u'}2^2+a^2\,\nrmS{u^{-1}}2^{-2}+\alpha\,\nrmS u2^2}{\nrmS up^2}\,.
\]
In what follows, we denote by $\mathrm H^1(\S^1)$ the subspace of the continuous functions $u$ on $(-\pi,\pi]$ such that $u(\pi)=u(-\pi)$ and $u'\in\mathrm L^2(\S^1)$. Notice that if $u\in\mathrm H^1(\S^1)$ is such that $u(s_0)=0$ for some $s_0\in(-\pi,\pi]$, then
\[
|u(s)|^2=\(\int_{s_0}^su'\,\mathrm ds\)^2\le\sqrt{2\pi}\,\nrmS{u'}2\,\sqrt{|s-s_0|}
\]
and $u^{-2}$ is not integrable. In this case we adopt the convention that $\mathcal Q_{a,p,\alpha}[u]=\mathsf Q_{p,\alpha}[u]$.
\begin{lem} For any $a\in(0,1/2)$, $p>2$, $\alpha>-\,a^2$,
\[
\mu_{a,p}(\alpha)=\min_{u\in\mathrm H^1(\S^1)\setminus\{0\}}\mathcal Q_{a,p,\alpha}[u]
\]
is achieved by a function $u>0$.\end{lem}
To prove this result, it is enough to check that the infimum~\eqref{minpsi} is achieved by a function $\psi\in X_a$ such that $\psi(s)\neq0$ for any $s\in(-\pi,\pi]$. Without loss of generality, we can assume that $\psi$ is an optimal function for~\eqref{Ineq:Interp} with $\nrmS\psi p=1$. Let us decompose $v(s)=\psi(s)\,e^{ias}$ as a real and an imaginary part, $v=v_1+i\,v_2$, which both solve the same Euler-Lagrange equation
\[
-\,v_j''+\alpha\,v_j=(v_1^2+v_2^2)^{\frac p2-1}\,v_j\,,\quad j=1\,,\;2\,.
\]
The Wronskian $w=(v_1\,v_2'-\,v_1'\,v_2)$ is constant.

Neither $v_1$ nor $v_2$ vanishes identically on~$\S^1$ because of~\eqref{Cdt:Per}. If both $v_1$ and $v_2$ vanish at the same point, then $w$ vanishes identically, which means that $v_1$ and $v_2$ are proportional. Again, this cannot be true because of the twisted boundary condition~\eqref{Cdt:Per}.\finprf

If $a=0$, $\mathcal Q_{a=0,p,\alpha}[u]=\mathsf Q_{p,\alpha}[u]$ for any $u\in\mathrm H^1(\S^1)\setminus\{0\}$.
\begin{lem}\label{lem:BE} For any $p>2$, if $0<\alpha\le1/(p-2)$, then $\mu_{0,p}(\alpha)=\alpha$ is achieved only by constant functions. Inequality~\eqref{Ineq:Interp} also holds with $p=-\,2$ and $\alpha=1/(p-2)=-1/4=\mu_{0,p}(-1/4)$, with equality achieved only by constant functions.\end{lem}
Both results (case $p>2$,~\cite{DoEsLa2013,2013arXiv1308.2259N} and case $p=-\,2$,~\cite{Exner-Harrell-Loss98}) were already known. For $p>2$, similar results have been found by various methods based on entropy techniques~\cite{MR2219352,MR2224869,DEKL2012,DEKL} and the \emph{carr\'e du champ} method of D.~Bakry and M.~Emery~\cite{Bakry-Emery85}. There are many earlier references~\cite{MR1784688,MR1230930,MR1231419,MR1132767,MR622225,MR0185064,zbMATH02502560} on \emph{Kolmogorov's inequalities}~\cite{MR0001787} corresponding to various other boundary conditions and intervals.

As a consequence of the cases $p>2$ and $p=-\,2$ we have the inequalities
\be{Ineq:InterpZero}
\nrmS{u'}2^2+\beta\,\nrmS u2^2\ge\beta\,\nrmS u p^2\quad\forall\,u\in\mathrm H^1(\S^1)
\ee
for any $p>2$ and $\beta\in(0,1/(p-2)]$, and
\be{Ineq:InterpZero2}
\nrmS{u'}2^2+\frac14\,\nrmS{{u}^{-1}}2^2\ge\frac14\,\nrmS u2^2\quad\forall\,u\in\mathrm H^1(\S^1)\,.
\ee
Inequality~\eqref{Ineq:InterpZero2} actually enters in the family of inequalities~\eqref{Ineq:InterpZero}, with the parameter $\beta=-1/4=1/(p-2)$ corresponding to the critical exponent $p=2\,d/(d-2)=-\,2$ since here $d=1$. This exponent is critical from the point of view of scalings because, at least for a function~$u$ with compact support in $(-\pi,\pi)$, $\nrmS up$ scales like $\nrmS{u'}2$. This is why a unified proof of both cases can be done with the Bakry-Emery method: see Appendix~\ref{Appendix:BE}.

\medskip\noindent We are now ready to study the key issues of Theorem~\ref{main}.
\begin{lem}\label{lem:rigidity} Let $p>2$, $a\in[0,1/2]$, and $\alpha>-\,a^2$.
\begin{enumerate}
\item[(i)] if $a^2\,(p+2)+\alpha\,(p-2)\le1$, then $\mu_{a,p}(\alpha)=a^2+\alpha$ and equality in~\eqref{Ineq:Interp} is achieved only by the constants,
\item[(ii)] if $a^2\,(p+2)+\alpha\,(p-2)>1$, then $\mu_{a,p}(\alpha)<a^2+\alpha$ and equality in~\eqref{Ineq:Interp} is not achieved by the constants.
\end{enumerate}
\end{lem}
In case (i), we can write
\begin{align*}
&\nrmS{u'}2^2+a^2\,\nrmS{u^{-1}}2^{-2}+\alpha\,\nrmS u2^2\\
&=(1-4\,a^2)\,\nrmS{u'}2^2+\alpha\,\nrmS u2^2\\
&\hspace*{12pt}+4\,a^2\(\nrmS{u'}2^2+\frac14\,\nrmS{{u}^{-1}}2^2\)
\end{align*}
and conclude using~\eqref{Ineq:InterpZero2} and then~\eqref{Ineq:InterpZero} with
\[
\beta=\frac{a^2+\alpha}{1-4\,a^2}\le\frac1{p-2}\,.
\]
In case (ii), let us consider the test function $u_\varepsilon:=1+\varepsilon\,w_1$, where $w_1$ is the eigenfunction corresponding to the first non zero eigenvalue of $-d^2/ds^2$ on $\mathrm H^1(\S^1)$, with Neumann boundary conditions, namely, $\lambda_1=1$ and \hbox{$w_1(s)=1+\cos s$}. A Taylor expansion shows that
\[
\mathcal Q_{a,p,\alpha}[u_\varepsilon]=a^2+\alpha+\big(1-a^2\,(p+2)-\,\alpha\,(p-2)\big)\,\varepsilon^2+o(\varepsilon^2)\,,
\]
which proves the result.\finprf
 
The proof of Lemma~\ref{lem:rigidity}, (i) relies on~\eqref{Ineq:InterpZero} and~\eqref{Ineq:InterpZero2}. It is remarkable that it does not use rigidity results based on the \emph{carr\'e du champ} method, at least directly. Notice that results similar to Lemma~\ref{lem:rigidity} were known for $p=+\infty$~\cite{MR1330606,Ilyin-Laptev-Loss-Zelik2016} using a Fourier representation of the operator and for an arbitrary $p>2$ if $\alpha=0$~\cite{2017arXiv171208829N}.

\medskip It follows from the definition of $\mathcal Q_{a,p,\alpha}[u]$ that $a\mapsto\mu_{a,p}(\alpha)$ is nondecreasing on $[0,1/2)$. The strict monotonicity follows from the existence of an optimal function, which is known by Lemma~\ref{Lem:Existence}. This concludes the proof of Theorem~\ref{main}. The remainder of this section is devoted to complementary results, which specify the range of $\mu_{a,p}(\alpha)$ when $a$ varies in $[0,1/2)$.

\medskip Let us consider
\[
\nu_p(\alpha):=\inf_{v\in\mathrm H_0^1(\S^1)\setminus\{0\}}\mathsf Q_{p,\alpha}[v]\,.
\]
Here $\mathrm H_0^1(\S^1)$ denotes the subspace of the functions $v\in\mathrm H^1(\S^1)$ such that $v(\pm\pi)=0$. Since~\eqref{Cdt:Per} is satisfied by any function in $\mathrm H_0^1(\S^1)$, we have the following estimate.
\begin{lem}\label{lem:Estim} If $p>2$, $\alpha>-\,a^2$ and $a\in\R$, then
\[
\mu_{a,p}(\alpha)\le\nu_p(\alpha)\,.
\]
Moreover, this inequality is strict if $a\in[0,1/2)$.\end{lem}

If $\{u_n\}_{n\in\N}$ is a minimizing sequence such that, for any $n\in\N$, $\nrmS{u_n}p=1$, then it is clearly bounded in $\mathrm H^1(\S^1)$, and so, by the compact Sobolev embeddings, it is relatively compact in $\mathrm L^2(\S^1)$, $\mathrm L^p(\S^1)$ and $C(\S^1)$. Up to subsequences, $\{u_n\}_{n\in\N}$ converges to some function $u$ weekly in $\mathrm H^1$ and strongly in $\mathrm L^2(\S^1)$, $\mathrm L^p(\S^1)$ and $C(\S^1)$. After noticing that $\mathsf Q_{p,\alpha}[\,|u|\,]=\mathsf Q_{p,\alpha}[u]$, we obtain the following result.
\begin{lem}\label{lem:Prop2.4} If $p>2$, $\alpha>-\,a^2$, then $\nu_p(\alpha)$ admits a non-negative minimizer.\end{lem}
The strict monotonicity of $a\mapsto\mu_{a,p}(\alpha)$ is a consequence of Lemma~\ref{lem:Prop2.4} and, as a consequence, we know that
\[
\mu_{a,p}(\alpha)<\mu_{1/2,p}(\alpha)\le\nu_p(\alpha)
\]
for any $a\in[0,1/2)$. It turns out that the last inequality is an equality.
\begin{thm}\label{main2} For any $p>2$ and $\alpha>-\,a^2$, we have
\[
\mu_{1/2,p}(\alpha)=\nu_p(\alpha)\,.
\]\end{thm}
This result was already known for the limit cases $p=2$,~\cite{Exner-Harrell-Loss98} and $p=+\infty$,~\cite{MR1330606,Ilyin-Laptev-Loss-Zelik2016}. To prove it, we set $v(s)=e^{is/2}\,\psi(s)$ and note that $v(s+2\pi)=-\,v(s)$ for all $s$, which follows from the periodicity condition~\eqref{Cdt:Per} with $a=1/2$. Moreover, the derivative $v'$ satisfies $v'(s+2\pi)=-\,v'(s)$. Note that these boundary conditions also hold for the real part and the imaginary part of $v$ separately. We call them $v_1$ and~$v_2$. Our problem is to minimize $\mathsf Q_{p,\alpha}[v]$ subject to these conditions. Both $v_1$ and $v_2$ must vanish at some point but \emph{a priori} these points need not be the same. We set $\eta_j=|v_j|$, $j=1$, $2$, and note that
\[
\mathsf Q_{p,\alpha}[v]=\frac{\int_{-\pi}^\pi\left[\eta_1'^2+\eta_2'^2\right]\,\mathrm d\sigma+\alpha\int_{-\pi}^\pi\left[\eta_1^2+\eta_2^2\right]\,\mathrm d\sigma}{\|\eta\|_p^2}\,.
\]
The functions $\eta_j$ are now periodic. They are not necessarily smooth but are at least continuous. Now we replace both $\eta_1$ and $\eta_2$ by their symmetric decreasing rearrangements around the point~$0$. The numerator decreases for the usual reasons and the denominator increases (see Lemma~\ref{lemmm1}, in Appendix~\ref{Appendix:symmetrization}). Thus, the symmetrically decreasing rearranged functions $\eta_1^*$ and $\eta_2^*$ have a maximum at $0$ and vanish at $\pm\pi$, so that $\eta_1^*+i\,\eta_2^*\in\mathrm H_0^1(\S^1)$. If $v$ is a minimizer of $\mathsf Q_{p,\alpha}$ under Condition~\eqref{Cdt:Per} with $a=1/2$, then
\[
\nu_p(\alpha)\le\mathsf Q_{p,\alpha}[\eta_1^*+i\,\eta_2^*]\le\mathsf Q_{p,\alpha}[v]=\mu_{1/2,p}(\alpha)\,.
\]
\finprf
With the convention that $\mathcal Q_{a,p,\alpha}[u]=\mathsf Q_{p,\alpha}[u]$ if $u\in\mathrm H_0^1(\S^1)$, we can claim that the infimum of $\mathcal Q_{a,p,\alpha}$ is attained by some $u\in\mathrm H^1(\S^1)\setminus\{0\}$ for any $a\in[0,1/2]$, including in the case $a=1/2$ for which the minimizer can be taken in $\mathrm H_0^1(\S^1)\setminus\{0\}$.

\section{Proof of Corollaries~\ref{Lambda1} and~\ref{Hardy} }\label{Proofcorollaries}

Let us start with the proof of Corollary~\ref{Lambda1}. Consider the quadratic form associated with $H_a-\varphi$. Using H\"older's inequality, we obtain
\begin{multline*}
\nrmS{\psi'+i\,a\,\psi}2^2-\int_{-\pi}^\pi\varphi\,|\psi|^2\,\mathrm d\sigma\\
\ge\nrmS{\psi'+i\,a\,\psi}2^2-\mu\,\nrmS\psi p^2
\end{multline*}
where $\mu=\nrmS\varphi q$ and $\frac1q+\frac2p=1$. Let us choose $\alpha$ such that $\mu_{a,p}(\alpha)=\mu$. It follows from~\eqref{Ineq:Interp} that
\[
\nrmS{\psi'+i\,a\,\psi}2^2-\mu\,\nrmS\psi p^2\ge-\,\alpha\,\nrmS\psi2^2
\]
and from Theorem~\ref{main} that $\mu_{a,p}(\alpha)=a^2+\alpha$ if $a^2\,(p+2)+\alpha\,(p-2)\le1$. This implies that
\[
\lambda_1(H_a-\varphi)\ge a^2-\nrmS\varphi q
\]
if $4\,a^2+\nrmS\varphi q\,(p-2)\le1$. In that case the equality is achieved by $\varphi\equiv const$. The proof is complete.\finprf

\medskip Now let us prove Corollary~\ref{Hardy}. Let $\mathbf x=(r,\vartheta)\in\R^2$ be polar coordinates in $\R^2$. Then we find
\begin{multline*}
\int_{\R^2}|(i\,\nabla+\mathbf a)\,\Psi|^2\,\mathrm d\mathbf x\\
=\int_0^\infty\int_{\S^1}\(r\,|\partial_r\Psi|^2+\frac1r\,|\partial_\vartheta\Psi+i\,a\,\Psi|^2\)\mathrm d\vartheta\,\mathrm dr\,.
\end{multline*}
Let $\tau>0$. Then
\begin{multline*}
\int_0^\infty\int_{\S^1}\frac1r\(|\partial_\vartheta\Psi+i\,a\,\Psi|^2-\tau\,\varphi\,|\Psi|^2\)\mathrm d\vartheta\,\mathrm dr\\
\ge\lambda_1\(H_a-\tau\,\varphi\)\int_0^\infty\int_{\S^1}\frac1r\,|\Psi|^2\,\mathrm d\vartheta\,\mathrm dr\\
\ge-\,\alpha_{a,p}(\tau\,\nrmS\varphi q)\int_0^\infty\int_{\S^1}\frac1r\,|\Psi|^2\,\mathrm d\vartheta\,.
\end{multline*}
Note that if $\tau=0$, then
\[
\alpha_{a,p}(\tau\,\nrmS\varphi q)=\alpha_{a,p}(0)=-\,a^2\,,
\]
and for a sufficiently large $\tau$ the value of $\alpha_{a,p}(\tau\,\nrmS\varphi q)$ is positive. Therefore we can find $\tau>0$ such that $\alpha_{a,p}(\tau\,\nrmS\varphi q)=0$. This value is unique since $\alpha_{a,p}(\mu)$ is strictly monotone with respect to $\mu$. The conclusion easily follows.\finprf

\appendix\section{A proof of Lemma~\ref{lem:BE} by the carr\'e du champ method}\label{Appendix:BE}

Let $\mathcal F_\beta[u]:=\nrmS{u'}2^2+\beta\,\big(\nrmS u2^2-\nrmS up^2\big)$. If $p>2$, it is enough to prove $\mu_{0,p}(\beta)=\beta$ for $\beta=\alpha_\star$, $\alpha_\star:=1/(p-2)$, because
\[
\mathcal F_\beta[u]=\big(1-\beta\,(p-2)\big)\,\nrmS{u'}2^2+\beta\,(p-2)\,\mathcal F_{\alpha_\star}[u]
\]
if $0<\beta\le\alpha_\star$. Let us consider a positive solution of the parabolic equation
\[
\frac{\partial u}{\partial t}=u''+(p-1)\,\frac{|u'|^2}u
\]
and compute
\begin{multline*}
-\frac d{dt}\mathcal F_{\alpha_\star}[u(t,\cdot)]=\intpi{\(|u''|^2-|u'|^2\)}\\
+\frac{p-1}3\intpi{\frac{|u'|^4}{u^2}}
\end{multline*}
using several integrations by parts. The first term in the r.h.s. is non-negative by the Poincar\'e inequality, as well as the second one. Notice that $\rho=|u|^p$ is a solution of the heat equation, so that positivity is preserved by the flow and $\mathcal F_{\alpha_\star}[u(t=0,\cdot)]\ge\lim_{t\to+\infty}\mathcal F_{\alpha_\star}[u(t,\cdot)]=0$, which is exactly~\eqref{Ineq:InterpZero} written with $u=u(t=0,\cdot)$. The strict positivity condition is easily removed by an approximation procedure. Exactly the same computations give the result in the case $p=-\,2$ and establish~\eqref{Ineq:InterpZero2}.

For $p>2$, the method is well known~\cite{Bakry-Emery85,DEKL}. The result for $p=-\,2$ was established earlier~\cite{Exner-Harrell-Loss98} but, as far as we know, this proof is new.

\section{A symmetrization result}\label{Appendix:symmetrization}

Here $f^*$ denotes the symmetric decreasing rearrangement of $f$.
\begin{lem}\label{lemmm1} Let $p\ge2$. For any non-negative functions $f$, $g\in\mathrm L^p(\S^1)$ we have that
\[
\int_{-\pi}^\pi\kern-6pt\(f^2+g^2\)^{p/2}\mathrm d\sigma\le\int_{-\pi}^\pi\({f^*}^2+{g^*}^2\)^{p/2}\mathrm d\sigma\,.
\]
\end{lem}
The case $p=2$ is obvious, in fact there is equality. Hence we assume that $p>2$. Write
\[
\(\int_{-\pi}^\pi\kern-6pt\(f^2+g^2\)^{p/2}\mathrm d\sigma\)^{2/p}=\sup_{\nrmS vq=1}\int_{-\pi}^\pi\kern-6pt\(f^2+g^2\)v\,\mathrm d\sigma
\]
where $1/q+2/p=1$. Clearly, we may choose $v$ to be positive. By standard rearrangement inequalities,
\[
\int_{-\pi}^\pi\kern-6pt\(f^2+g^2\)v\,\mathrm d\sigma\le\int_{-\pi}^\pi\({f^*}^2+{g^*}^2\)v^*\,\mathrm d\sigma
\]
and $\nrmS{v^*}q=\nrmS vq$: the proof is completed with
\begin{multline*}
\sup_{\nrmS{v^*}q=1}\int_{-\pi}^\pi\({f^*}^2+{g^*}^2\)v^*\,\mathrm d\sigma\\
=\(\int_{-\pi}^\pi\({f^*}^2+{g^*}^2\)^{p/2}\mathrm d\sigma\)^{2/p}\,.
\end{multline*}

\section{Some numerical results}\label{Appendix:numerics}

\begin{figure}[ht]
\includegraphics[width=8cm]{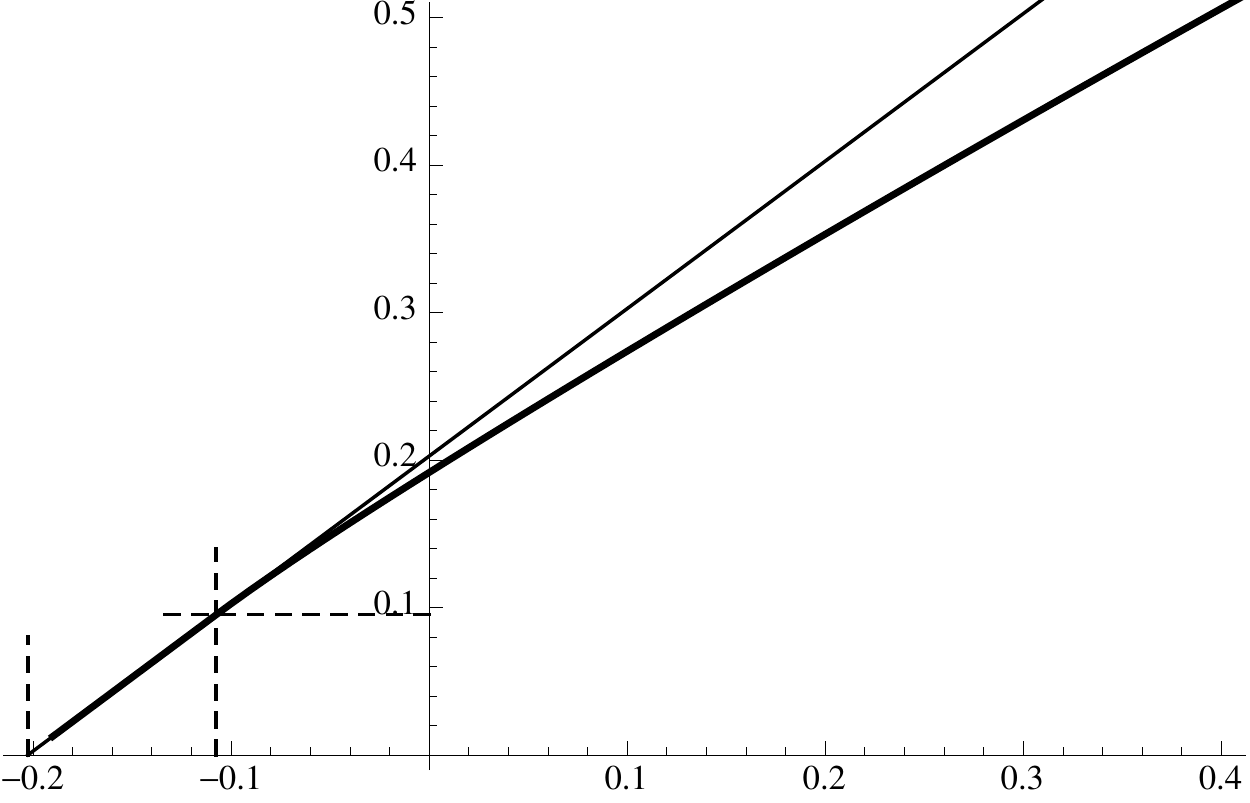}
\caption{\label{DoEsLaLo-F1} The curve $\alpha\mapsto\mu_{a,p}(\alpha)$ with $p=4$ and $a=0.45$. The only solutions to~\eqref{EL} are the constant functions for any $\alpha$ such that $-\,a^2=-\,0.2025\le\alpha\le-\,0.1075$ and, in this range, $\mu_{a,p}(\alpha)=a^2+\alpha$. A branch of non-constant optimizers of~\eqref{Ineq:Interp} bifurcates at $\alpha=-\,0.1075$.}
\end{figure}

\begin{figure}[ht]
\includegraphics[width=8cm]{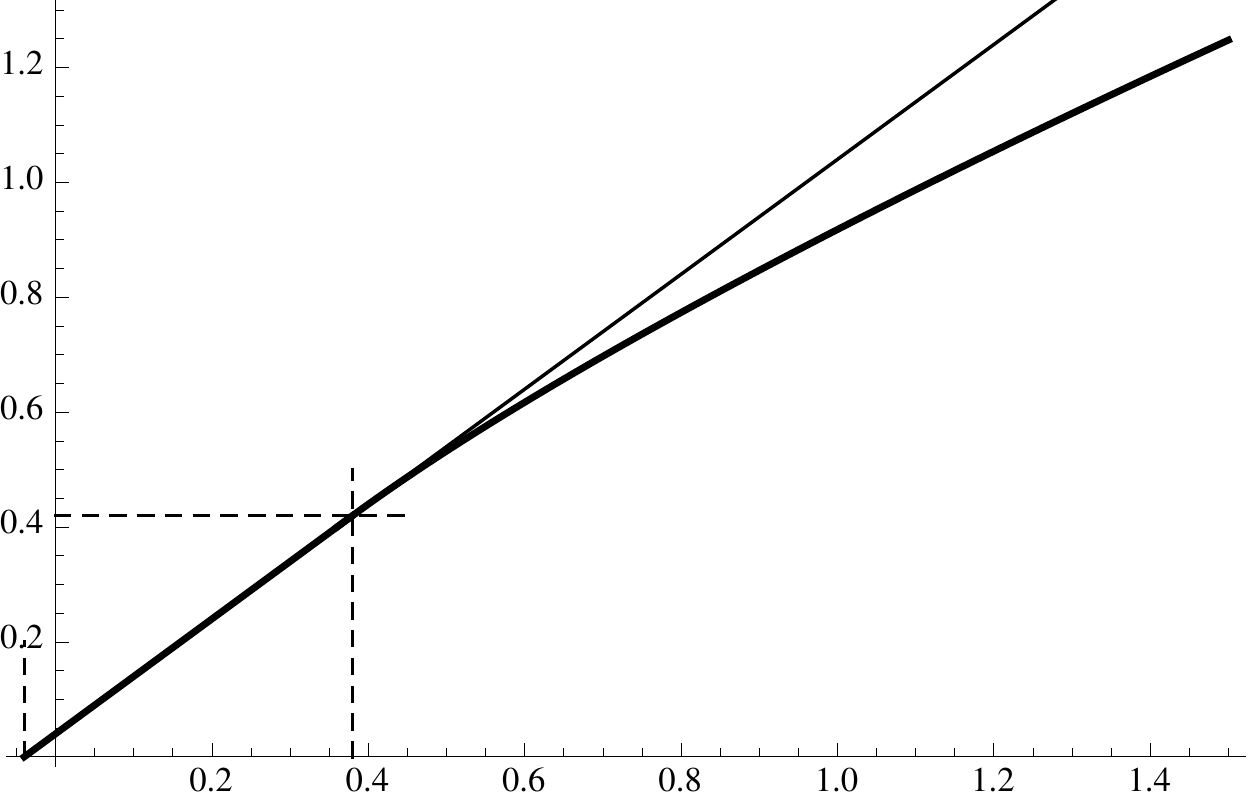}
\caption{\label{DoEsLaLo-F2} The curve $\alpha\mapsto\mu_{a,p}(\alpha)$ with $p=4$ and $a=0.2$. Here the branch of non-constant optimizers of~\eqref{Ineq:Interp} bifurcates at $\alpha=0.38$ which corresponds to $a^2\,(p+2)+\alpha\,(p-2)=1$.}
\end{figure}

To compute the curve $\alpha\mapsto\mu_{a,p}(\alpha)$, we systematically solve the Euler-Lagrange equation associated with the variation of $\mathcal Q_{a,p,\alpha}$, \emph{i.e.},
\be{EL}
-\,u''+\frac{a^2}{u^3\(\intpi{\frac1{u^2}}\)^2}+\alpha\,u=u^{p-1}
\ee
where the solution $u>0$ is normalized by
\[
\mu_{a,p}(\alpha)\(\intpi{u^p}\)^{\frac2p-1}=1\,.
\]
This condition \emph{a posteriori} provides the numerical value of $\mu_{a,p}(\alpha)$. To impose the boundary conditions $u'(0)=u'(\pi)=0$, we use a shooting method and solve~\eqref{EL} on~$\R$ with the conditions $u'(0)=0$ and $u(0)=\lambda>0$. To emphasize the dependence in $\lambda$, let us denote it by $u_\lambda$. For any $\lambda>0$, $\lambda\neq(a^2+\alpha)^{1/(p-2)}$, the solution is non-constant and periodic so that
\[
\rho(\lambda)=\min\{s>0\,:\,u_\lambda'(s)=0\}
\]
is well defined. The shooting parameter $\lambda$ is then determined by the condition that $\rho(\lambda)=\pi$. Since~\eqref{EL} involves a nonlocal term, an additional fixed-point procedure is needed to adjust the coefficient of $u^{-3}$ in the equation. Some plots are shown in Figs.~\ref{DoEsLaLo-F1} and~\ref{DoEsLaLo-F2}.

\begin{figure}[ht]
\includegraphics[width=8cm]{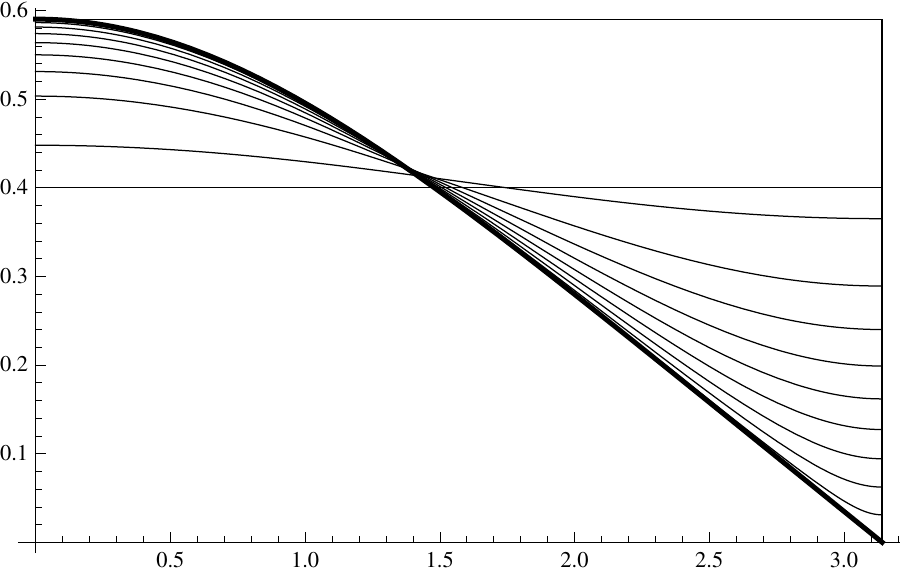}
\caption{\label{DoEsLaLo-F3} Here $p=4$ and $\alpha=0$. Plot of the solution of~\eqref{EL} for $a=0.40$, $0.41$,\ldots $0.49$. The thick curve solves $u''+u^{p-1}=0$ and it is explicit. Similar patterns are found when $\alpha\neq0$, with a non-explicit curve solving~\eqref{ELlimit} in the limit as $a\to1/2$.}
\end{figure}

Equality in~\eqref{Ineq:Interp} is achieved only by constant functions according to Lemma~\ref{lem:rigidity} if $a^2\,(p+2)+\alpha\,(p-2)\le1$: in this case, $\lambda=(a^2+\alpha)^{1/(p-2)}\equiv u_\lambda$. For any $a\in(0,1/2)$ such that $a^2\,(p+2)+\alpha\,(p-2)>1$, our method provides us with a non-constant solution $u$ of~\eqref{EL} which realizes the equality in~\eqref{Ineq:Interp}. As $a\to1/2$, the integral $\intpi{u^{-2}}$ diverges, so that the limit curve is described by the solution of 
\be{ELlimit}
-\,u''+\alpha\,u=u^{p-1}
\ee
with boundary conditions $u'(0)=0$ and $u(\pi)=0$. See Fig.~\ref{DoEsLaLo-F3}.

AL was partially supported by RSF grant No. 18-11-00032 
\begin{acknowledgments}\begin{spacing}{0.8}{\noindent\small This research has been partially supported by the project \emph{EFI}, contract~ANR-17-CE40-0030 (J.D.) of the French National Research Agency (ANR), by RSF grant No. 18-11-00032 (A.L.) and by the NSF grant DMS-DMS-1600560 (M.L.). Some of the preliminary investigations were done at the Institute Mittag-Leffler during the fall program \emph{Interactions between Partial Differential Equations \& Functional Inequalities}. The authors thank A. I. Nazarov for pointing them several important references.\\
\copyright\,2018 by the authors. This paper may be reproduced, in its entirety, for non-commercial purposes.}\end{spacing}\end{acknowledgments}
\section*{References}
\bibliography{DoEsLaLo}
\end{document}